\documentclass[english,11pt]{article}
\usepackage[latin1]{inputenc}
\usepackage{amsmath,amsthm,amssymb,babel}

\newtheorem{teo}{Theorem}[section]
\newtheorem{defin}{Definition}[section]
\newtheorem{prop}{Proposition}[section]

\newtheorem{lemma}{Lemma}[section]

\def\proof{{\it Proof.}\ }

\def\eq#1{(\ref{#1})}
\def\neweq#1{\begin{equation}\label{#1}}
\def\endeq{\end{equation}}

\def\RR{\mathbb R}
\def\NN{\mathbb N}
\def\phi{\varphi}

\def\di{\displaystyle}
\def\ri{\rightarrow}

\def\r2{{\mathbb R}^{2}}

\def\eq#1{(\ref{#1})}

\def\phi{\varphi}
\def\RR{\mathbb R}
\def\NN{\mathbb N}

\def\di{\displaystyle}
\def\ri{\rightarrow}
\def\intom{\int_\Omega}
\def\huo{H^1_0(\Omega )}
\def\incep{\left\{\begin{array}{cl} }
 \def\termin{\end{array}\right. }
\def\2af{2^*_\alpha}

\def\proof{{\it Proof.}\ }

\title{\sc Subcritical perturbations of resonant linear problems 
 with sign-changing potential}
\author{\sc Teodora-Liliana Dinu\\
\small Department of Mathematics,
``Fra\c tii Buze\c sti" College, 200352 Craiova, Romania\\ \small Email : {\tt
tldinu@gmail.com}}

\date{}
\begin{document}
\maketitle

\begin{abstract}
 We establish  existence and multiplicity theorems for a
 Dirichlet boundary value problem at resonance, which is a nonlinear subcritical
 perturbation of a linear eigenvalue problem studied by Cuesta. Our framework includes 
 a sign-changing potential and we locate the solutions by using the Mountain Pass
lemma and the Saddle Point theorem.
\end{abstract}
\normalsize

\section{Introduction and main results}

Let $\Omega$ be an arbitrary open set in ${\mathbb R}^{N}$, $N\geq 2$, and assume that %%@
$V:\Omega\ri\RR$ is a variable potential.
Consider
 the eigenvalue problem
\neweq{eig}
 \left\{\begin{tabular}{ll}
&$-\Delta u=\lambda V(x)u$ \quad  $\mbox{in}\
\Omega\,,$\\
&$u\in\huo$. \\
\end{tabular} \right.
\end{equation}
Problems of this type have a long history. If $\Omega$ is bounded
and $V\equiv 1$, problem \eq{eig} is related to the Riesz-Fredholm theory
of self-adjoint and compact operators (see, e.g., Theorem VI.11 in \cite{6}). The case of a %%@
non-constant
potential $V$ has been first considered in the pioneering papers
of Bocher \cite{bo}, Hess and Kato \cite{hk}, Minakshisundaran and  Pleijel \cite{mp} and Pleijel %%@
\cite{p}.
For instance, Minakshisundaran and  Pleijel \cite{mp}, \cite{p} studied
the case where $\Omega$ is bounded, $V\in L^\infty (\Omega)$,
$V\geq 0$ in $\Omega$
and $V>0$ in $\Omega_0\subset\Omega$ with $|\Omega_0|>0$. An important
contribution in the study of Problem \eq{eig}
if $\Omega$ and $V$ are not necessarily bounded has been given recently
by Cuesta \cite{cuesta} (see also
Szulkin and Willem \cite{sw})
under the assumption that the sign-changing potential $V$
satisfies
$$
V^+\not=0\ \mbox{ and }\ V\in L^s(\Omega)\,,\leqno (H)$$
where $s>N/2$ if $N\geq 2$ and $s=1$ if $N=1$.
As usually, we have denoted   $V^{+}(x)=\max\{V(x),0\}$. Obviously,
$V=V^+-V^-$, where $V^{-}(x)=\max\{-V(x),0\}$.

In order to study the main properties (isolation, simplicity) of the principal eigenvalue of %%@
\eq{eig},
Cuesta \cite{cuesta} proved that
the minimization problem
$$\min\left\{\intom |\nabla u|^2dx;\ u\in \huo ,\ \intom V(x)u^2dx=1\right\}$$
has a positive solution $\varphi_1=\varphi_1 (\Omega)$ which is an eigenfunction of \eq{eig}
corresponding to the eigenvalue $\lambda_1:=\lambda_1(\Omega)=
\int_\Omega |\nabla\varphi_1 |^2dx$.

Our purpose in this paper is to study the existence of
solutions of the perturbed nonlinear boundary value problem
\begin{eqnarray}\label{et1}
\left\{\begin{array}{ll}
&\di -\Delta u=\lambda_1V(x)u+g(x,u) \quad \mbox{in  } \Omega, \\
 &\di u=0 \quad \mbox{on } \partial \Omega,\\
 &\di u\not\equiv 0\quad \mbox{in  } \Omega,
\end{array}\right.
\end{eqnarray}
where $V$ satisfies $(H)$ and
$g:\Omega\times{\mathbb R}\rightarrow{\mathbb R}$ is
a Carath\'eodory function satisfying $g(x,0)=0$ and with subcritical growth, that is,
\begin{eqnarray}
 \di    |g(x,s)|\le a_{0}\cdot |s|^{r-1} +b_{0}, \qquad \mbox{for all }
 s\!\in{\mathbb R},\ \ \mbox{a.e. }
     \, x\!\in\Omega,
\end{eqnarray}
 for some constants $a_{0}$, $b_{0}>0$, where $1\le r<2^*$. We recall that $2^*$ denotes the %%@
critical Sobolev exponent, that is, $2^*:=\frac{2N}{N-2}$ if $N\geq 3$ and $2^*=+\infty$ if %%@
$N\in\{1,2\}$.

Problem \eq{et1} is resonant at infinity and equations of this type
have been first studied by Landesman and Lazer \cite{ll} in connection
with concrete problems arising in Mechanics.

Set $G(x,s)=\int\limits_{0}^{s}g(x,t)dt$.
Throughout this paper we assume that
 there exist $k,\ m\in L^{1}(\Omega)$, with $m\geq 0$, such that
\begin{eqnarray}\label{et5} |G(x,s)|\le k(x),\qquad\mbox{for all } s\in{\mathbb R},\ \ %%@
\mbox{a.e.}\ x\in\Omega\,;
\end{eqnarray}
\begin{eqnarray}\label{et7}
\liminf\limits_{s\ri 0}\frac{G(x,s)}{s^{2}}=m(x),\qquad\mbox{a.e.}\ x\in\Omega\,.
\end{eqnarray}

The energy functional associated to Problem \eq{et1} is
$$F(u)=\frac 12\intom\left(|\nabla u|^2-\lambda_1V(x)u^2\right)dx-\intom G(x,u)dx\,,$$
for all $u\in H^1_0(\Omega)$.

From the variational characterization of
 $\lambda_{1}$ and using (\ref{et5}) we obtain
$$F(u)\ge-\int_{\Omega}G(x,u(x))dx\ge-|k|_{1}>-\infty\,,$$
for all $u\in H^1_0(\Omega)$
and, consequently, $F$ is bounded from below. Let us consider $u_{n}=\alpha_{n}\varphi_{1}$,
where $\alpha_n\ri\infty$.
Then
the estimate $\int_{\Omega}|\nabla\varphi_{1}|^{2}dx=\lambda_{1}\int_{\Omega}\varphi_{1}^{2}dx$
yields $F(u_{n})=-\int_{\Omega}G(x,\alpha_{n}\varphi_{1})dx\le|k|_{1}<\infty$. Thus,
$\lim_{n\ri\infty}F(u_{n})<\infty$.
Hence the sequence $(u_{n})_{n}\subset H^1_0(\Omega)$ defined by $ u_{n}=\alpha_{n}\varphi_{1}$
satisfies
 $|\!|u_{n}|\!|\ri\infty$ and $F(u_{n})$ is bounded.
In conclusion, if we suppose that (\ref{et5})
holds true then the energy functional $F$ is bounded from below and is not coercive.

Our first result is the following.

\begin{teo}\label{theo1}
Assume that for all $\omega\subset\Omega$  with  $|\Omega\setminus\omega|>0$ we have
\begin{equation}\label{et8}
\int\limits_{\omega}\limsup\limits_{|s|\ri\infty}G(x,s)dx\le 0\quad\mbox{  and  }\quad
\int\limits_{\Omega\setminus\omega}G(x,s)dx\le 0
\end{equation}
and
\begin{equation}\label{et9}
\int\limits_{\Omega}\limsup\limits_{|s|\ri\infty}G(x,s)dx\le 0\,.
\end{equation}
Then Problem (\ref{et1}) has at lest one solution.
\end{teo}

Denote $V:=\mbox{Sp}\,(\varphi_{1})$. Since $1={\rm dim}\,V<\!\infty$, there exists
a closed complementary subspace
$W$ of $V$, that is, 
$W\cap V=\{0\}$ and
$H^1_0(\Omega)=V\oplus W$.
For such a closed complementary subspace $W\!\subset\! H^1_0(\Omega)$,
 denote
$$\lambda_{W}:=\inf\left\{\frac{\int_{\Omega}|\nabla w|^{2}dx}{\int_{\Omega}w^{2}dx};\ 
                              w\in W,\ 
                              w\not=0
                              \right\}\,.$$

  The following result establishes a multiplicity result, provided
  $G$ satisfies a certain subquadratic condition. 
							  
\begin{teo}\label{theo2}
Assume that the conditions of Theorem \ref{theo1} are fulfilled and the following additional %%@
assumption is fulfilled:  
\begin{eqnarray}\label{et10}
G(x,s)\le\frac{\lambda_{W}-\lambda_{1}}{2}\,s^{2},\qquad\mbox{for all } s\in{\mathbb R},\ a.e.\ %%@
x\in\Omega\,.
\end{eqnarray}
Then Problem \eq{et1} has at least two solutions.
\end{teo}

In the next two theorems, we prove the existence of a solution under
the following assumptions on the potential $G$:
$$\limsup\limits_{|s|\ri\infty}\frac{G(x,s)}{|s|^{q}}
\le b<\infty \quad\hbox{uniformly a.e. }\ x\!\in\!\Omega\,, \ q>2;\leqno (G_{1})_{q}$$
$$\liminf\limits_{|s|\ri\infty}
\frac{g(x,s)s-2G(x,s)}{|s|^{\mu}}\ge a>0\quad\hbox{uniformly a.e. }\ x\!\in\!\Omega;\leqno %%@
(G_{2}^{+})_{\mu}$$
$$\limsup\limits_{|s|\ri\infty}
\frac{g(x,s)s-2G(x,s)}{|s|^{\mu}}\le -a<0\quad\hbox{uniformly a.e. }\ x\!\in\!\Omega\,.\leqno %%@
(G_{2}^{-})_{\mu}$$

\begin{teo}\label{theo3}
Assume that $G$ satisfies conditions $(G_{1})_{q}$, $(G_{2}^{+})_{\mu}$ [or $(G_{2}^{-})_{\mu}$] %%@
and
$$\limsup\limits_{s\ri 0}\frac{2G(x,s)}{s^{2}}\le\alpha<\lambda_{1}<\beta\le\liminf
\limits_{|s|\ri\infty}\frac{2G(x,s)}{s^{2}}\quad\mbox{uniformly a.e. }\ x\in\Omega\,,\leqno %%@
(G_{3})$$
with $\mu>2N/(q-2)$ if $N\ge 3$ or $\mu>q-2$ if $1\le N\le 2$.
Then Problem \eq{et1} has at least one solution.
\end{teo}

\begin{teo}\label{theo4}
Assume that $G(x,s)$ satisfies $(G_{2}^{-})_{\mu}$ [or $(G_{2}^{+})_{\mu}$],
for some $\mu>0$, and
$$\lim\limits_{|s|\ri\infty}\frac{G(x,s)}{s^{2}}=0\quad\mbox{uniformly a.e. }
\ x\in\Omega\,.\leqno (G_{4})$$
Then Problem \eq{et1} has at least one solution.
\end{teo}

\section{Compactness conditions and auxiliary results}

Let $E$ be a  reflexive real Banach space with norm $|\!|\cdot|\!|$ and let
 $I:E\rightarrow {\mathbb R}$ be a $C^{1}$ functional. 
We assume that there exists
a compact embedding $E\hookrightarrow X$, where $X$ is a real Banach space,
and that the following interpolation type inequality holds:
$$ |\!|u|\!|_{X}\le\psi (u)^{1-t}|\!|u|\!|^{t}\,, \qquad\mbox{for all } u\in E\,,\leqno (H_{1})$$
for some $t\in (0,1)$ and some homogeneous function $\psi :E\ri{\mathbb R}_{+}$
of degree one. An example of such a framework is the following:
$E=H^1_0(\Omega)$, $X=L^{q}(\Omega)$, $\psi (u)=|u|_{\mu}$,
where $0<\mu <q<2^{*}$. Then, by the interpolation inequality (see \cite[Remarque 2, p.~57]{6}) we %%@
have 
  $$|u|_{q}\le |u|_{\mu}^{1-t}|u|_{2^{*}}^{t}\,, \quad
    \mbox{where  } \frac{1}{q} =\frac{1-t}{\mu}+\frac{t}{2^{*}}\,. $$
The Sobolev inequality yields
$|u|_{2^{*}}\le c|\!|u|\!|$, for all $u\in H^1_0(\Omega)$. Hence
 $$|u|_{q}\le k|u|_{\mu}^{1-t}|\!|u|\!|^{t}\,, \qquad\mbox{for all } u\in H^1_0(\Omega)$$
and this is a $(H_{1})$ type inequality.
 
 We recall below the following Cerami compactness conditions.

\begin{defin}\label{defcerami}
a) The functional $I:E\ri {\mathbb R}$ is said to satisfy condition $(C)$ at the level %%@
$c\in{\mathbb R}$ [denoted $(C)_{c}$]
if any sequence $(u_{n})_{n}\subset E$ such that $I(u_{n})\ri c$
and $(1+|\!|u_{n}|\!|)\cdot |\!|I^{'}(u_{n})|\!|_{E^{*}}\ri 0$ possesses a
convergent subsequence.\\
b) The functional $I:E\ri {\mathbb R}$ is said to satisfy condition $(\hat C)$ at the level
$c\in{\mathbb R}$ [denoted $(\hat C)_{c}$] if any sequence $(u_{n})_{n}\subset E$ such that
$I(u_{n})\ri c$ and $(1+|\!|u_{n}|\!|)\cdot|\!|I^{'}(u_{n})|\!|_{E^{*}}\ri 0$
possesses a bounded subsequence.
\end{defin}

We observe that the above conditions are weaker than the usual Palais-Smale condition
$(PS)_{c}$: any sequence $(u_{n})_{n}\subset E$ such that
$I(u_{n})\ri\! c$ and $|\!|I^{'}(u_{n})|\!|_{E^{*}}\ri\! 0$ possesses a convergent
subsequence.

 Suppose that $I(u)=J(u)-N(u)$, where $J$ is $2$-homogeneous
 and $N$ is not $2$-homogeneous at infinity. We recall that $J$ is $2$-homogeneous if
$J(\tau u)=\tau^{2} J(u)$, for all $\tau\in{\mathbb R}$ and for any $u\in E$.
We also recall that the functional $N\in C^{1}(E,{\mathbb R})$ is said to be not $2$-homogeneous
at infinity if there exist $a$, $c>0$ and  $\mu >0$ such that
 $$|\langle N^{'}(u),u\rangle-2N(u)|\ge a\psi (u)^{\mu}-c\,, \qquad\mbox{for all } 
 u\in E\,.\leqno (H_{2})$$

We introduce the following additional hypotheses on the functionals $J$ and $N$:
$$ J(u)\ge k|\!|u|\!|^{2}\,, \qquad\mbox{for all } u\in E\leqno (H_{3})$$
$$|N(u)|\le b|\!|u|\!|_{X}^{q}+d\,, \qquad\mbox{for all } u\in E\,,\leqno (H_{4})$$
for some constants $k$, $b$, $d>0$ and $q>2$.

\begin{teo}\label{auxi}
Assume that assumptions $(H_{1})$, $(H_{2})$, $(H_{3})$ and $(H_{4})$ are fulfilled, with $qt<2$. %%@
Then the
functional $I$ satisfies condition $(\hat C)_{c}$, for all $c\in{\mathbb R}$.
\end{teo}

\proof Let $(u_{n})_{n}\subset E$ such that $I(u_{n})\ri c$ and
$(1+|\!|u_{n}|\!|)|\!|I^{'}(u_{n})|\!|_{E^{*}}\ri 0$. We have
$$
\begin{array}{ll}
\di |\langle I^{'}(u),u\rangle-pI(u)|&\di=|\langle J^{'}(u)-N^{'}(u),u\rangle-2J(u)+pN(u)|\\
&\di =|\langle J^{'}(u),u\rangle-2J(u)-(\langle N^{'}(u),u\rangle-2N(u))|\,.
\end{array}
$$
But $J$ is $2$-homogeneous and
$$ \frac{J(u+tu)-J(u)}{t}=J(u)\,\frac{(1+t)^{2}-1}{t}\,.$$
This implies $\langle J^{'}(u),u\rangle=2J(u)$ and $$|\langle I^{'}(u),u\rangle-2I(u)|=|\langle %%@
N^{'}(u),u\rangle-2N(u)|\,.$$
From $(H_{2})$ we obtain
$$ |\langle I^{'}(u),u\rangle-2I(u)|=|\langle N^{'}(u),u\rangle-2N(u)|\ge a\psi (u)^{\mu}-c\,.$$
Letting $u=u_{n}$ in the inequality from above we have:
 $$a\psi(u_{n})^{\mu}\le c+|\!|I^{'}(u_{n})|\!|_{E^{*}}|\!|u_{n}|\!|+2|I(u_{n})|\,.$$
Thus, by our hypotheses, for some $c_{0}>0$ and all positive integer $n$,
$\psi(u_{n})\le c_{0}$ and hence, the sequence $\{\psi(u_n)\}$ is bounded.
 Now, from $(H_{1})$ and $(H_{4})$ we obtain
$$ J(u_{n})=I(u_{n})+N(u_{n})\le b|\!|u_{n}|\!|_{X}^{q}+d_{0}\le b\psi(u_{n})^{(1-t)q}
|\!|u_{n}|\!|^{qt}+d_{0}\,.$$
Hence
$$ J(u_{n})\le b_{0}|\!|u_{n}|\!|^{qt}+d_{0}\,, \qquad\mbox{for all } n\in\NN\,,$$
for some $b_{0}$, $d_{0}>0$. Finally, $(H_{3})$ implies
$$c|\!|u_{n}|\!|^{2}\le b_{0}|\!|u_{n}|\!|^{qt}+d_{0}\,,\qquad\mbox{for all } n\in\NN\,.$$
Since $qt<2$, we conclude that $(u_{n})_{n}$ is bounded in $E$.\qed

\begin{prop}\label{eqqq}
Assume that $I(u)=J(u)-N(u)$ is as above, where  $N^{'}:E\ri E^{*}$ is a compact
operator and $J^{'}:E\ri E^{*}$ is an isomorphism from $E$ onto $J^{'}(E)$.
Then  conditions $(C)_{c}$ and $(\hat C)_{c}$ are equivalent.
\end{prop}

\proof It is enough to show that $(\hat C)_{c}$ implies $(C)_{c}$.
Let $(u_{n})_{n}\subset E$ be a sequence such that $I(u_{n})\ri c$ and
$(1+|\!|u_{n}|\!|)|\!|I^{'}(u_{n})|\!|_{E^{*}}\ri 0$. From $(\hat C)_{c}$ we
obtain  a bounded subsequence $(u_{n_{k}})_{k}$ of $(u_{n})_{n}$. But $N^{'}$
is a compact operator.
Then $N^{'}(u_{n_{k_{l}}})\stackrel{l}{\ri}f^{'}\in E^{*}$,
where $(u_{n_{k_{l}}})$ is a subsequence of $(u_{n_{k}})$. Since
$ (u_{n_{k_{l}}})$ is a bounded sequence and
$(1+|\!|u_{n_{k_{l}}}|\!|)|\!|I^{'}(u_{n_{k_{l}}})|\!|_{E^{*}}\ri 0$, it
follows that $ |\!|I^{'}(u_{n_{k_{l}}})|\!|\ri 0$. Next, using the relation
$$\di u_{n_{k_{l}}}=J^{'^{-1}}(I^{'}(u_{n_{k_{l}}})+N^{'}(u_{n_{k_{l}}}))\,,$$
we obtain that $(u_{n_{k_{l}}})$ is a convergent subsequence of $(u_{n})_{n}$.\qed

\section{Proof of Theorem \ref{theo1}}

We first show that the energy functional
$F$ satisfies the Palais-Smale condition at level $c<0$:
any sequence $(u_{n})_{n}\subset H^1_0(\Omega)$ such that
$F(u_{n})\ri c$ and
$|\!|F^{'}(u_{n})|\!|_{H^{-1}}\ri 0$ possesses a convergent subsequence.

Indeed, it suffices to show that such a sequence $(u_{n})_{n}$
has a bounded subsequence (see the Appendix). Arguing by contradiction, we suppose that %%@
$|\!|u_{n}|\!|\ri\infty$. We distinguish the following two distinct situations.

\smallskip
{\sc Case 1}: $|u_{n}(x)|\ri\infty$ a.e. $x\in\Omega$.
Thus, by our hypotheses,
$$\begin{array}{ll} c&\di=\liminf\limits_{n\ri\infty}F(u_{n})=\liminf\limits_{n\ri\infty}
\left\{\frac{1}{2}\int\limits_{\Omega}|\nabla u_{n}|^{2}dx-\frac{\lambda_{1}}{2}
\int\limits_{\Omega}u_{n}^{2}dx-\int\limits_{\Omega}G(x,u_{n}(x))dx\right\}\\
&\di\ge\liminf\limits_{n\ri\infty}\left(-\int\limits_{\Omega}G(x,u_{n}(x)))dx\right)=
-\limsup\limits_{n\ri\infty}\int\limits_{\Omega}G(x,u_{n}(x))dx\\
&\di=-\limsup\limits_{|s|\ri\infty}\int\limits_{\Omega}G(x,s)dx.\end{array}$$
Using Fatou's lemma we obtain
$$\limsup\limits_{|s|\ri\infty}\int\limits_{\Omega}G(x,s)dx\le
\int\limits_{\Omega}\limsup\limits_{|s|\ri\infty}G(x,s)dx\,.$$
Our assumption (\ref{et9}) implies $c\ge 0$. This is a contradiction
because $c<0$. Therefore $(u_{n})_{n}$ is bounded in $H^1_0(\Omega)$.

\smallskip
{\sc Case 2}: there exists $\omega\subset\subset\Omega$ such that $|\Omega\setminus\omega|>0$
and $|u_{n}(x)|\not\ri\infty\ \ \mbox{for all } x\in\Omega\setminus\omega$.
It follows that there exists a subsequence, still denoted by $(u_{n})_{n}$, which is bounded
in $\Omega\setminus\omega$. So, there exists $k>0$
such that $|u_{n}(x)|\le k$, for all $x\in\Omega\setminus\omega$. Since $I(u_{n})\ri c$ we obtain %%@
there exists $M$ such that $I(u_{n})\le M$, for all $n$. We have
$$\frac{1}{2}\,|\!|u_{n}|\!|^{2}-\frac{\lambda_{1}}{2}\,|u_{n}|_{2}^{2}-|k|_{1}\le
I(u_{n})\le M\qquad\mbox{ as }\quad|\!|u_{n}|\!|\ri\infty\,.$$
It follows that $|u_{n}|_{2}\ri\infty$. We have 
$$\di|u_{n}|_{2}^{2}=\int\limits_{\Omega}u_{n}^2dx=
\int\limits_{\Omega\setminus\omega}u_{n}^{2}dx+\int\limits_{\omega}u_{n}^{2}dx
\le k^{2}|\Omega\setminus\omega|+|u_{n}|_{L^{2}(\omega)}^{2}\,.$$
This shows that $|u_{n}|_{L^{2}(\omega)}\ri\infty$.
If $(u_{n})_{n}$ is bounded in $\omega$, this implies $(|u_{n}|_{L^{2}(\omega)})_{n}$ is
bounded, and this is a contradiction. Therefore $|u_{n}(x)|\ri\infty$ for all $x\in\omega$.
So, by Fatou's lemma and our assumption (\ref{et8}),
$$
\begin{array}{ll}
c&\di=\liminf\limits_{n\ri\infty}F(u_{n})\ge
-\limsup\limits_{n\ri\infty}\int\limits_{\Omega}G(x,u_{n}(x))dx\\
&\di=-\limsup\limits_{n\ri\infty}\left(\int\limits_{\Omega\setminus\omega}G(x,u_{n}(x))dx+
\int\limits_{\omega}G(x,u_{n}(x))dx\right)\\
&\di\ge-\limsup\limits_{n\ri\infty}\int\limits_{\Omega\setminus\omega}G(x,u_{n}(x))dx-
\limsup\limits_{n\ri\infty}\int\limits_{\omega}G(x,u_{n}(x))dx\\
&\di\ge-\limsup\limits_{n\ri\infty}\int\limits_{\Omega\setminus\omega}G(x,u_{n}(x))dx-
\int\limits_{\omega}\limsup\limits_{|s|\ri\infty}G(x,s)dx\ge 0\,.
\end{array}
$$
This implies $c\ge 0$  which contradicts our hypothesis $c<0$. This contradiction
shows that $(u_{n})_{n}$ is bounded in $H^1_0(\Omega)$, and hence
$F$ satisfies the Palais-Smale condition at level $c<0$.

The assumption (\ref{et7}) is equivalent with: there exist $\delta_{n}\searrow 0$ and %%@
$\varepsilon_{n}\in L^{1}(\Omega)$ with 
$|\varepsilon_{n}|_{1}\ri 0$ such that
\begin{equation}\label{et12}
\int\limits_{\Omega}\frac{G(x,s)}{s^{2}}dx\ge\int\limits_{\Omega}m(x)dx-
\int\limits_{\Omega}\varepsilon_{n}(x)dx\,,\qquad\mbox{for all } 0<|s|\le\delta_{n}\,.
\end{equation}
But $|\varepsilon_{n}|_{1}\ri 0$ implies
that for all $\varepsilon>0$ there exists $n_{\varepsilon}$ such that for all 
$n\ge n_{\varepsilon}$ we have $|\varepsilon_{n}|_{1}<\varepsilon$.
Set $\varepsilon=\int_{\Omega}m(x)\varphi_{1}^{2}dx/|\!|\varphi_{1}|\!|_{L^\infty}^{2}$
and fix $n$  large enough so that
$$L:=\int\limits_{\Omega}m(x)\varphi_{1}^{2}(x)dx-
|\varepsilon_{n}|_{1}|\!|\varphi_{1}|\!|_{L^\infty}^{2}>0\,.$$
Take $v\in V$ such that $|\!|v|\!|\le\delta_{n}/|\!|\varphi_{1}|\!|_{L^{\infty}}$.
We have $F(v)=-\int\limits_{\Omega}G(x,v(x))dx$.
The inequality (\ref{et12}) is equivalent with
$$\int\limits_{\Omega}G(x,s)dx\ge\int\limits_{\Omega}m(x)s^{2}dx-\int\limits_{\Omega}\varepsilon_{%%@
n}(x)s^{2}dx$$
and therefore
\begin{eqnarray}\label{et13}
F(v)\!=\!-\int\limits_{\Omega}G(x,v(x))dx\!\le\!-\int\limits_{\Omega}m(x)v^{2}(x)dx+
\int\limits_{\Omega}\varepsilon_{n}(x)v^{2}(x)dx\,.
\end{eqnarray}
By our choice of $v$ in $V=\mbox{Sp}\,(\varphi_{1})$ we have
$$|v(x)|=|\alpha|\,|\varphi_1(x)|\le|\alpha||\!|\varphi_{1}|\!|_{L^{\infty}}\le
|\alpha|\frac{\delta_{n}}{|\!|v|\!|}\,.$$
But, from (\ref{et13}),
$$
\begin{array}{ll}
F(v)&\di\le-\int\limits_{\Omega}mv^{2}dx+\int\limits_{\Omega}\varepsilon_{n}v^{2}dx\le
-\int\limits_{\Omega}m|\alpha|^{2}\varphi_1^{2}dx+|\alpha|^{2}\int\limits_{\Omega}\varepsilon_{n}|%%@
\!|\varphi_1|\!|_{L^{\infty}}^{2}dx\\
&\di=|\alpha|^{2}\left(-\int\limits_{\Omega}m\varphi_1^{2}dx+
|\varepsilon_{n}|_{1}|\!|\varphi_1|\!|_{L^{\infty}}^{2}\right)=
-L|\alpha|^{2}=-L|\!|v|\!|^{2}\,,\end{array}$$
with $n$ so large that $L>0$.
Therefore we obtain the existence of some $v_{0}\in V$ such that $F(v_{0})<0$.
This implies $l=\inf_{H^1_0(\Omega)}F<0$. But the functional $F$
satisfies the Palais-Smale condition (P-S)$_{c}$, for all $c<0$. This implies that there exists
$u_{0}\in H^1_0(\Omega)$ such that $F(u_{0})=l$. In conclusion, $u_{0}$
is a critical point of $F$ and consequently it is a  solution of Problem (\ref{et1}).
Our assumption $g(x,0)=0$ implies $F(0)=0$ and we know that $F(u_{0})=l<0$, that is,
$u_{0}\not\equiv 0$. Therefore $u_{0}\in H^1_0(\Omega)$ is a nontrivial
solution of (\ref{et1}) and the proof of Theorem \ref{theo1} is complete.\qed

\section{Proof of Theorem \ref{theo2}}
Let $X$ be a real Banach space and $F:X\ri{\mathbb R}$ be a $C^{1}$-functional. Denote
$$K_{c}:=\{u\in X;\ F^{'}(u)=0\mbox{ and }F(u)=c\}$$
$$F^{c}:=\{u\in X;\ F(u)\le c\}\,.$$
The proof of Theorem \ref{theo2} makes use of the following deformation lemma (see \cite{11}).

\begin{lemma}
\label{deformm}
Suppose that $F$ has no critical values in the interval $(a,b)$ and that
$F^{-1}(\{a\})$ contains at most a finite number of critical points of $F$.
Assume that the Palais-Smale condition $(P-S)_{c}$ holds for every $c\in[a,b)$.
Then there exists an
$F$-decreasing homotopy of homeomorphism $h:[0,1]\times F^{b}\setminus K_{b}\ri X$
such that
$$h(0,u)=u\,,\qquad\mbox{for all } u\in F^{b}\setminus K_{b}$$
$$h(1,F^{b}\setminus K_{b})\subset F^{a}$$ and
$$h(t,u)=u\,,\qquad\mbox{for all } u\in F^{a}\,.$$
\end{lemma}

We are now in position to give the proof of Theorem \ref{theo2}.
Fix $n$  large enough so that
$$F(v)\le-L|\!|v|\!|^{2}\,,\qquad\mbox{for all } v\in V\mbox{ with %%@
}|\!|v|\!|\le\frac{\delta_{n}}{|\!|\varphi_1|\!|_{L^{\infty}}}\,.$$
Denote  $d:=\sup_{\partial B}F$, where $B=\{v\in V;\ |\!|v|\!|\le R\}$ and
$R=\delta_{n}/|\!|\varphi_1|\!|_{L^{\infty}}$.
We suppose that 0 and $u_{0}$ are the only critical points of $F$ and we
show that this yields a contradiction. For any $w\in W$ we have
$$F(w)=\frac{1}{2}\left(\int\limits_{\Omega}|\nabla %%@
w|^{2}dx-\lambda_{1}\int\limits_{\Omega}w^{2}dx\right)
-\int\limits_{\Omega}G(x,w(x))dx\,.$$
Integrating in (\ref{et10}) we find
\begin{eqnarray}\label{et14}
-\int\limits_{\Omega}G(x,w(x))dx\ge\frac{\lambda_{1}-\lambda_{W}}{2}\int\limits_{\Omega}w^{2}dx\,.
\end{eqnarray}
Combining the definition of $\lambda_{W}$ with relation (\ref{et14}) we obtain
\begin{eqnarray}\label{et15}
\begin{array}{ll}
\di F(w)&\di\ge\frac{1}{2}\int\limits_{\Omega}|\nabla %%@
w|^{2}dx-\frac{\lambda_{1}}{2}\int\limits_{\Omega}w^{2}dx+
\frac{\lambda_{1}-\lambda_{W}}{2}\int\limits_{\Omega}w^{2}dx\\
&\di=\frac{1}{2}\left(\int\limits_{\Omega}|\nabla %%@
w|^{2}dx-\lambda_{W}\int\limits_{\Omega}w^{2}dx\right)\ge 0\,.
\end{array}
\end{eqnarray}
Using $0\in W$, $F(0)=0$ and relation (\ref{et15}) we find $\inf_{W}F=0$.
If $v\in\partial B$ then $F(v)\le-LR<0$ and, consequently,
$$d=\sup\limits_{\partial B}F<\inf\limits_{W}F=0\,.$$
Obviously, $$l=\inf\limits_{H^1_0(\Omega)}F\le\inf\limits_{\partial B}F<d=\sup\limits_{\partial %%@
B}F\,.$$
Denote 
$$\alpha :=\inf\limits_{\gamma\in\Gamma}\sup\limits_{u\in B}F(\gamma(u))\,,$$
where $$\Gamma :=\{\gamma\in C(B,H^1_0(\Omega));\ \gamma(v)=v\quad\mbox{for all } v\in\partial %%@
B\}\,.$$
It is known (see the Appendix) that $\gamma(B)\cap W\not=\emptyset$, for all $\gamma\in\Gamma$.
Since $\inf\limits_{W}F=0$, we have $F(w)\ge 0$ for all $w\in W$.
Let $u\in B$ such that $\gamma(u)\in W$. It follows that $F(\gamma(u))\ge 0$ and hence
 $\alpha\ge 0$. The Palais-Smale condition holds true at level $c<0$ and
the functional $F$ has no critical value in the interval $(l,0)$, So, by
 Lemma \ref{deformm}, we obtain an $F$ decreasing homotopy $h:[0,1]\times F^{0}\setminus K_{0}\ri %%@
H^1_0(\Omega)$
such that
$$
\begin{array}{lll}
\di h(0,u)=u\,,\qquad\mbox{for all } u\in F^{0}\setminus K_{0}=F^{0}\setminus\{0\}\,;\\
\di h(1,F^{0})\setminus\{0\}\subset F^{l}=\{u_{0}\}\,;\\
\di h(t,u)=u\,,\qquad\mbox{for all } u\in F^{l}\,.
\end{array}
$$
Consider the application $\gamma_{0}:B\ri H^1_0(\Omega)$ defined by
\begin{eqnarray*}
\gamma_{0}=\left\{\begin{array}{ll}
         \di u_{0}\,,\qquad\mbox{if } |\!|v|\!|<R/2\\
         \di h\left(\frac{2(R-|\!|v|\!|)}{R},\frac{Rv}{2|\!|v|\!|}\right)\,,\qquad\mbox{if %%@
}|\!|v|\!|\ge R/2\,.
           \end{array}\right.
\end{eqnarray*}
Since $\gamma_{0}(v)=h(1,v)=u_{0}$ if  $|\!|v|\!|=R/2$, we deduce that
$\gamma_{0}$ is continuous.

If $v\in\partial B$ then $v=R\varphi_1$ and $F(R\varphi_1)\le 0$. Then
$v\in F^{0}\setminus\{0\}$ and hence $\gamma_{0}(v)=v$. Therefore we obtain
that $\gamma_{0}\in\Gamma$. The condition that $h$ is $F$ decreasing is equivalent
with
$$s>t\qquad\mbox{implies  }F(h(s,u))<F(h(t,u))\,.$$
Let us consider $v\in B$. We distinguish the following situations.

\smallskip
{\sc Case 1}: $|\!|v|\!|<\frac{R}{2}$.
In this case, $\gamma_{0}(v)=u_{0}$ and $F(u_{0})=l<d$.

\smallskip
{\sc Case 2}: $|\!|v|\!|\ge\frac{R}{2}$.
If $|\!|v|\!|=R/2$ then $\gamma_{0}(v)=h(1,v)$ and if $|\!|v|\!|=R$
then $\gamma_{0}(v)=h(0,v)$.
But $0\le t\le 1$ and $h$ is $F$ decreasing. It follows that
$$F(h(0,v))\ge F(h(t,v))\ge F(h(1,v))\,,$$
that is, $F(\gamma_{0}(v))\le F(h(0,v))=F(v)\le d$.

From these two cases we obtain $F(\gamma_{0}(v))\le d$, for all $v\in B$
and from the definition of $\alpha$ we have $0\le\alpha\le d<0$. This is a
contradiction. We conclude that  $F$ has a another critical
point $u_{1}\in H^1_0(\Omega)$ and, consequently, Problem (\ref{et1}) has a
second nontrivial weak solution. \qed

\section{Proof of Theorems \ref{theo3} and \ref{theo4}}

We will use the following classical critical point theorems.

\begin{teo}
{\bf(Mountain Pass, \cite{1})}. Let $E$ be a real Banach space.
Suppose that $I\in C^{1}(E,{\mathbb R})$ satisfies condition
$(C)_{c}$, for all $c\in{\mathbb R}$ and, for some $\rho>0$ and $u_{1}\in E$ with
$|\!|u_{1}|\!|>\rho$,
$$\max\{I(0),I(u_{1})\}\le\hat\alpha<\hat\beta\le\inf\limits_{|\!|u|\!|=\rho}I(u)\,.$$
Then $I$ has a critical value $\hat c\ge\hat\beta$,
characterized by
$$\hat c=\inf\limits_{\gamma\in\Gamma}\max\limits_{0\le\tau\le1}I(\gamma(\tau))\,,$$
where $\Gamma :=\{\gamma\in C([0,1],E);\ \gamma(0)=0,\gamma(1)=u_{1}\}$.
\end{teo}

\begin{teo}
{\bf(Saddle Point, \cite{9})}. Let $E$ be a real Banach space.
Suppose that $I\!\in\! C^{1}(E,{\mathbb R})$ satisfies
condition $(C)_{c}$, for all  $c\in{\mathbb R}$ and, for some $R>0$ and some
$E=V\oplus W$ with ${\rm dim}\,V<\infty$,
$$ \max\limits_{v\in V,|\!|v|\!|=R}I(v)\le\hat\alpha<\hat\beta\le
\inf\limits_{w\in W}I(w)\,.$$
Then $I$ has a critical value $\hat c\ge\hat\beta$,
characterized by
$$\hat c=\inf\limits_{h\in\Gamma}\max\limits_{v\in V,|\!|v|\!|\le R}I(h(v))\,,$$
where $\Gamma=\{h\in C(V\bigcap\bar B_{R},E);\ h(v)=v,\,\mbox{for all }
 v\in\partial B_{R}\}$.
\end{teo}

\begin{lemma}\label{lema41}
Assume that $G$ satisfies conditions $(G_{1})_{q}$ and $(G_{2}^{+})_{\mu}$ [or %%@
$(G_{2}^{-})_{\mu}$], with
 $\mu\!>2N/(q-2)$ if $N\geq 3$ or $\mu\!>\!q-2$ if $1\!\le\! N\!\le 2$.
Then the functional $F$ satisfies condition $(C)_{c}$ for all $c\in{\mathbb R}$.
\end{lemma}

\proof Let $$\di %%@
N(u)=\frac{\lambda_{1}}{2}\int\limits_{\Omega}u^{2}dx+\int\limits_{\Omega}G(x,u)dx \qquad
\mbox{and}\qquad
J(u)=\frac{1}{2}|\!|u|\!|^{2}\,.$$ 
Obviously, $J$ is homogeneous of degree $2$
and $J^{'}$ is an isomorphism of $E=H^1_0(\Omega)$ onto $ J^{'}(E)\subset
H^{-1}(\Omega)$.
It is known that $N^{'}:E\ri E^{*}$ is a compact operator.
Proposition \ref{eqqq} ensures that conditions $(C)_{c}$ and
$(\hat C)_{c}$ are equivalent. So, it suffices to show that $(\hat C)_{c}$ holds
for all $c\in{\mathbb R}$. Hypothesis $(H_{3})$ is trivially satisfied, 
whereas $(H_{4})$ holds
true from $(G_{1})_{q}$. Condition $(G_{1})_{q}$ implies that
$$\inf\limits_{|s|>0}\sup\limits_{|t|>|s|}\frac{G(x,t)}{|t|^{q}}\le b\,.$$
Therefore there exists $s_{0}\not=0$ such that
$$\sup\limits_{|t|>|s_{0}|}\frac{G(x,t)}{|t|^{q}}\le b\qquad\mbox{and}\qquad
G(x,t)\le b|t|^{q},\quad\mbox{for all } t\mbox{ with }|t|>|s_{0}|\,.$$
The boundedness is provided by the continuity of the application
$[-s_{0},s_{0}]\ni t\longmapsto G(x,t)$. It follows that $\int\limits_{\Omega}G(x,u)dx\le %%@
b|u|_{q}^{q}+d$.
But $\ N(u)=\lambda_{1}|u|_{2}^{2}/2+\int_{\Omega}G(x,u)dx$
and $q>2$. If $|u|_{q}\le 1$ then we obtain $(H_{4})$; if not, we have
$|u|_{2}\le k|u|_{q}$ because $\Omega$ is bounded.
Therefore $|u|_{2}^{2}\le k|u|_{q}^{2}\le k|u|_{q}^{q}$ and
finally $(H_{4})$ is fulfilled. Hypothesis $(H_{1})$ is a direct 
consequence of the Sobolev
inequality. It remains to show that
hypothesis $(H_{2})$ holds true, that is, the functional $N$ is not $2$-homogeneous
at infinity. Indeed, using assumption $(G_{2}^{+})_{\mu}$ (a similar argument works if %%@
$(G_{2}^{-})_{\mu}$ is fulfilled) together with the
subcritical condition on $g$ yields
$$\sup\limits_{|s|>0}\inf\limits_{|t|>|s|}\frac{g(x,t)t-2G(x,t)}{|t|^{\mu}}\ge a>0\,.$$
It follows that there exists $s_{0}\not=0$ such that
$$ \inf\limits_{|t|>|s_{0}|}
\frac{g(x,t)t-2G(x,t)}{|t|^{\mu}}\ge a\,.$$
Hence
$$g(x,t)t-2G(x,t)\ge a|t|^{\mu}\,,\qquad\mbox{for all }|t|>|s_{0}|\,.$$
The application $t\mapsto g(x,t)t-2G(x,t)$ is continuous in $[-s_{0},s_{0}]$,
therefore it is bounded.
We obtain $g(x,t)-2G(x,t)\ge a_{1}|t|^{\mu}-c_{1}$, for all  $s\in{\mathbb R}$ and a.e. %%@
$x\in\Omega$.
We deduce that
$$
\begin{array}{ll}
\di |\langle N^{'}(u),u\rangle-2N(u)|&\di=\left|\int\limits_{\Omega}(g(x,u)u-2G(x,u))dx\right|\\
&\di \ge a_{1}|\!|u|\!|_{\mu}^{\mu}-c_{2}\,,
\qquad\mbox{for all } u\in H^1_0(\Omega)\,.\end{array}$$
Consequently, the functional $N$ is not $2$-homogeneous at infinity.

Finally, when $N\ge 3$, we observe that condition $\mu>N(q-2)/2$ is
equivalent with $\mu>2^{*}(q-2)/{2^{*}-2}$.
From $1/q=(1-t)/\mu+t/2^{*}$ we obtain
$(1-t)/\mu=(2^{*}-qt)/(2^{*}q)$. Hence
$(2^{*}-qt)/q<(1-t)(2^{*}-2)/(q-2)$ and, consequently,
$(q-2^{*})(2-tq)<0$. But $q<2^{*}$ and this implies $2>tq$.
Similarly, when $1\le N\le 2$, we choose some $2^{**}>2$ sufficiently large so that
$\mu>2^{**}(q-2)/(2^{**}-2)$ and $t\in(0,1)$ be as above.
The proof of Lemma is complete in view of  Theorem \ref{auxi}.\qed

Our next step is to show that condition $(G_{3})$ implies the geometry of
the Mountain Pass theorem for the functional $F$.

\begin{lemma}\label{dinnoulema}
Assume that $G$ satisfies the hypotheses
$$
\limsup\limits_{|s|\ri\infty}\frac{G(x,s)}{|s|^{q}}\le b<\infty
\quad\mbox{uniformly }a.e.\ x\!\in\!\Omega\leqno (G_{1})_{q}$$
$$\limsup\limits_{s\ri 0}\frac{2G(x,s)}{s^{2}}\le\alpha<\lambda_{1}<\beta\le\liminf
\limits_{|s|\ri\infty}\frac{2G(x,s)}{|s|^{2}}\qquad\mbox{uniformly a.e. }\ x\in\Omega\,.\leqno %%@
(G_{3})$$
Then there exists $\rho$, $\gamma>0$ such that $F(u)\ge\gamma$ if $|u|=\rho$.
Moreover, there exists $\varphi_{1}\in H^1_0(\Omega)$ such that
$F(t\varphi_{1})\ri -\infty$ as $t\ri\infty$.
\end{lemma}

\proof In view of our hypotheses and the subcritical growth condition, we
obtain
$$\liminf\limits_{|s|\ri\infty}\frac{2G(x,s)}{s^{2}}\ge\beta\mbox{  is equivalent with  }
\sup\limits_{s\not=0}\inf\limits_{|t|>|s|}\frac{2G(x,t)}{t^{2}}\ge\beta\,.$$
There exists $s_{0}\not=0$ such that $\inf\limits_{|t|>|s_{0}|}\frac{2G(x,t)}{t^{2}}\ge\beta$
and therefore $\frac{2G(x,t)}{t^{2}}\ge\beta$, for all $|t|>|s_{0}|$  or
 $G(x,t)\ge\frac{1}{2}\beta t^{2}$, provided $|t|>|s_{0}|$.
We choose $t_{0}$  such that $|t_{0}|\le|s_{0}|$ and $G(x,t_{0})<\frac{1}{2}\beta|t_{0}|^{2}$.
Fix $\varepsilon>0$. There exists $B(\varepsilon,t_{0})$ such that
$G(x,t_{0})\ge\frac{1}{2}(\beta-\varepsilon)|t_{0}|^{2}-B(\varepsilon,t_{0})$.
Denote $B(\varepsilon)=\sup_{|t_{0}|\le|s_{0}|}B(\varepsilon,t_{0})$.
 We obtain for any given
$\varepsilon>0$ there exists $B=B(\varepsilon)$ such that
\begin{eqnarray}
G(x,s)\ge\frac{1}{2}\,(\beta-\varepsilon)|s|^{2}-B\,,\qquad\mbox{for all } s\in{\mathbb R}\,,\ %%@
\mbox{a.e. }\ x\in\Omega\,.
\end{eqnarray}
Fix arbitrarily $\varepsilon>0$. In the same way, using the second inequality of
$(G_{3})$ and $(G_{1})_{q}$ it follows that there exists $A=A(\varepsilon)>0$
such that
\begin{eqnarray}\label{et17}
2G(x,t)\le(\alpha+\varepsilon)t^{2}+2(b+A(\varepsilon))|t|^{q}\,,\qquad
\mbox{for all } t\in{\mathbb R}\,,\ \mbox{a.e. }\ x\in\Omega\,.
\end{eqnarray}
We now choose $\varepsilon>0$ so that $\alpha+\varepsilon<\lambda_{1}$ and we
use  (\ref{et17}) together with the Poincar\'e inequality
  to obtain the first assertion of the lemma.
  
Set  $H(x,s)=\lambda_{1}s^{2}/2+G(x,s)$. Then $H$  satisfies
$$
\limsup\limits_{|s|\ri\infty}\frac{H(x,s)}{|s|^{q}}\le b<\infty\,,\qquad
\mbox{uniformly a.e. } \ x\in\Omega\leqno (H_{1})_{q}$$
$$\limsup\limits_{s\ri 0}\frac{2H(x,s)}{s^{2}}\le\alpha<\lambda_{1}
<\beta\le\liminf\limits_{|s|\ri\infty}\frac{2H(x,s)}{s^{2}}\,,\qquad \mbox{uniformly a.e. }
\ x\!\in\!\Omega\,.\leqno (H_{3})$$
In the same way, for any given $\varepsilon>0$ there exists $A=A(\varepsilon)>0$
and $B=B(\varepsilon)$ such that
\begin{eqnarray}\label{et18}
\begin{array}{ll}
\di\frac{1}{2}(\beta-\varepsilon)s^{2}-B&\di\le H(x,s)\\
&\di\le\frac{1}{2}(\alpha+\varepsilon)s^{2}+A|s|^{q}\,,\qquad\mbox{for all } s\in{\mathbb R}\,,\ %%@
\mbox{a.e.}\ x\in\Omega\,.
\end{array}
\end{eqnarray}
We have
$$\begin{array}{ll}
\di F(u)&\di=\frac{1}{2}|\!|u|\!|^{2}-\int\limits_{\Omega}H(x,u)dx\ge\frac{1}{2}|\!|u|\!|^{2}
-\frac{1}{2}(\alpha+\varepsilon)|u|_{2}^{2}-A|u|_{q}^{q}\\
&\di\ge\frac{1}{2}\left(1-\frac{\varepsilon+\alpha}{\lambda_{1}}\right)
|\!|u|\!|^{2}-Ak|\!|u|\!|^{q}\,.\end{array}$$
We can assume without loss of generality that $q>2$. Thus, the above estimate 
yields $F(u)\ge\gamma$ for some $\gamma>0$, as long as $\rho>0$ is small, thus
proving the first assertion of the lemma.

On the other hand, choosing now $\varepsilon>0$ so that $\beta-\varepsilon>\lambda_{1}$
and using (\ref{et18}), we obtain
$$F(u)\le\frac{1}{2}|\!|u|\!|^{2}-\frac{\beta-\varepsilon}{2}\,|u|_{2}^{2}+B|\Omega|\,.$$
We consider $\varphi_{1}$ be the $\lambda_{1}$-eigenfunction with $|\!|\varphi_{1}|\!|=1$.
It follows that
$$F(t\varphi_{1})\le\frac{1}{2}\left(1-\frac{\beta-\varepsilon}{\lambda_{1}}\right)t^{2}
+B|\Omega|\ri -\infty\qquad \mbox{as }t\ri\infty.$$
This proves the second assertion of our lemma.\qed

\begin{lemma}\label{betty}
Assume that $G(x,s)$ satisfies the conditions $(G_{2}^{-})_{\mu}$ (for some $\mu>0$) and
$$
\lim\limits_{|s|\ri\infty}\frac{G(x,s)}{s^{2}}=0\,,\qquad
\mbox{uniformly a.e. }\ x\in\Omega\,.\leqno (G_{4})$$
 Then there exists a subspace $W$ of $H^1_0(\Omega)$
such that  $H^1_0(\Omega)=V\oplus W$ and\\
(i) $F(v)\ri -\infty$, as $|\!|v|\!|\ri\infty$, $v\in V$;\\
(ii) $F(w)\ri\infty$, as $|\!|w|\!|\ri\infty$, $w\in W$.
\end{lemma}

\proof
(i) The condition $(G_{2}^{-})_{\mu}$ is equivalent with
$$\exists\  s_{0}\not=0\mbox{ such that }g(x,s)s-2G(x,s)\le -a|s|^{\mu}\,,\quad\mbox{for all }
|s|\ge|s_{0}|=\! R_{1}\,,\ \mbox{a.e.}\ x\in\Omega\,.$$
Integrating the identity
$$\frac{d}{ds}\frac{G(x,s)}{|s|^{2}}=\frac{g(x,s)s^{2}-2|s|G(x,s)}
{s^{4}}=\frac{g(x,s)|s|-2G(x,s)}{|s|^{3}}$$
over an interval $[t,T]\subset[R,\infty)$ and using the above inequality we find
$$\frac{G(x,T)}{T^{2}}-\frac{G(x,t)}{t^{2}}\le -a\int\limits_{t}^{T}s^{\mu-3}ds=
\frac{a}{2-\mu}\left(\frac{1}{T^{2-\mu}}-\frac{1}{t^{2-\mu}}\right)\,.$$
Since we can assume that $\mu<2$ and using the above relation, we obtain
$$G(x,t)\ge\hat a t^{\mu}\,,\qquad\mbox{for all } t\ge\! R_{1}\,,\  \mbox{where }\hat %%@
a=\frac{a}{2-\mu}>0\,.$$
Similarly, we show that
$$G(x,t)\ge\hat a|t|^{\mu}\,,\qquad\mbox{for }|t|\ge R_{1}\,.$$
Consequently, $\lim_{|t|\ri\infty}G(x,t)\!=\!\infty$. Now, letting $v=t\varphi_{1}\in V$
and using the variational characterization of $\lambda_{1}$, we have
$$F(v)\ge-\int\limits_{\Omega}G(x,v)dx\ri-\infty\,,\qquad\mbox{as  }
|\!|v|\!|=|t||\!|\varphi_{1}|\!|\ri\infty\,.$$
This result is a consequence of the Lebesgue's dominated convergence theorem.

(ii) Let $V=$Sp$\,(\varphi_{1})$ and $W\!\subset\! H^1_0(\Omega)$ be a closed
complementary subspace to $V$. Since $\lambda_{1}$ is an
eigenvalue of Problem \eq{eig}, it follows that there exists $d>0$
such that
$$\inf_{0\not=w\in W}
\frac{\int\limits_{\Omega}|\nabla\!w|^{2}dx}{\int\limits_{\Omega}|w|^{2}dx}
\ge\lambda_{1}+d\,.$$
Therefore
$$|\!|w|\!|^{2}\ge(\lambda_{1}+d)|w|_{2}^{2}\,,\qquad\mbox{for all } w\in W\,.$$
Let $0\!<\!\varepsilon\!<\!d$. From $(G_{4})$ we deduce that
there exists $\delta\!=\!\delta(\varepsilon)>0$  such that  for all $s$ satisfying $|s|\!>\delta$
we have
$2G(x,s)/s^{2}\le\varepsilon$, a.e. $x\in\Omega$.
In conclusion 
$$G(x,s)-\frac{1}{2}\varepsilon s^{2}\le M\,,\qquad\mbox{for all } s\in{\mathbb R}\,,$$
where  $$M:=\sup_{|s|\le\delta}\left(G(x,s)-\frac{1}{2}\,\varepsilon\, s^{2}\right)<\infty\,.$$
Therefore
\begin{eqnarray*}
\begin{array}{ll}
\di F(w)&\di=\frac{1}{2}|\!|w|\!|^{2}-
\frac{\lambda_{1}}{2}|w|_{2}^{2}-\int\limits_{\Omega}G(x,w)dx\\ 
&\di\ge\frac{1}{2}|\!|w|\!|^{2}-\frac{\lambda_{1}}{2}|w|_{2}^{2}-
\frac{1}{2}\varepsilon|w|_{2}^{2}-M\\ 
&\di\ge\frac{1}{2}\left(1-\frac{\lambda_{1}+\varepsilon}{\lambda_{1}+d}\right)
|\!|w|\!|^{2}-M=N|\!|w|\!|^{2}
-M\,,\qquad\mbox{for all } w\in W\,.
\end{array}
\end{eqnarray*}
It follows that $F(w)\ri\infty$ as $|\!|w|\!|\ri\infty$, for all $w\in W$,
which completes the proof of the lemma. \qed

\medskip
{\sc Proof of Theorem \ref{theo3}}. In view of Lemmas \ref{lema41} and \ref{dinnoulema}, we may %%@
apply the Mountain Pass theorem
with $u_{1}=t_{1}\varphi_{1}$, $t_{1}>0$ being such that $F(t_{1}\varphi_{1})\le 0$ (this is
possible from Lemma \ref{dinnoulema}). Since $F(u)\ge\gamma$ if $|\!|u|\!|=\rho$, we have %%@
$$\max\{F(0),F(u_{1})\}=0=\hat\alpha<\inf\limits_{|\!|u|\!|=\rho}
F(u)=\hat\beta\,.$$
It follows that the energy functional $F$ has a critical value $\hat c\ge\hat\beta>0$
and, hence, Problem (\ref{et1}) has a nontrivial solution $u\in H^1_0(\Omega)$.
\qed

\medskip
{\sc Proof of Theorem \ref{theo4}}.
In view of Lemmas \ref{lema41} and \ref{betty}, we may apply the Saddle Point theorem with
$\hat\beta:=\inf_{w\in W}F(w)$ and $R>0$ being such that
$\sup_{|\!|v|\!|=R}F(v):=\hat\alpha<\hat\beta$, for all $v\in V$ (this is possible
because $F(v)\ri-\infty$ as $|\!|v|\!|\ri\infty$).
It follows that  $F$ has a critical value $\hat c\ge\hat\beta$, which is a weak solution of
Problem (\ref{et1}). \qed

\section{Appendix}
Throughout this section we assume that $\Omega\subset{\mathbb R}^{N}$ is a bounded domain with %%@
smooth boundary.
We start with the following auxiliary result.

\begin{lemma}\label{brezislema}
Let 
$g:\Omega\times{\mathbb R}\ri{\mathbb R}$ be a Carath\'eodory function and assume that
there exist some constants $a$, $b\ge 0$ such that
$$|g(x,t)|\le a+b|t|^{r/s}\,,\qquad\mbox{for all } t\in{\mathbb R}\,,\ \mbox{a.e.}\ %%@
x\in\Omega\,.$$
Then the application $\varphi(x)\mapsto g(x,\varphi(x))$ is in
$C(L^{r}(\Omega),L^{s}(\Omega))$.
\end{lemma}

\proof
For any $u\in L^{r}(\Omega)$  we have
$$\int\limits_{\Omega}|g(x,u(x))|^{s}dx\le\int\limits_{\Omega}(a+b|u|^{r/s})^{s}dx\le
 2^{s}\int\limits_{\Omega}(a^{s}+b^{s}|u|^{r})dx\le c\int\limits_{\Omega}(1+|u|^{r})dx<\infty
 \,.$$
This shows that if $\varphi\in L^{r}(\Omega)$ then $g(x,\varphi)\in L^{s}(\Omega)$.
Let $u_{n}$, $u\in L^{r}$ be such that $|u_{n}-u|_{r}\ri 0$.
By Theorem IV.9 in \cite{6}, there exist a subsequence $(u_{n_{k}})_{k}$ and
$h\in L^{r}$ such that
$\di u_{n_{k}}\ri u$ a.e. in $\Omega$ and
$|u_{n_{k}}|\le h$ a.e. in $\Omega$.
By our hypotheses it follows that
$\di g(u_{n_{k}})\ri g(u)$ a.e. in $\Omega$. Next, we observe that
$$ |g(u_{n_{k}})|\le a+b|u_{n_{k}}|^{r/s}\le
a+b|h|^{r/s}\in L^{s}(\Omega)\,.$$ 
So, by Lebesgue's dominated convergence theorem,
$$\di |g(u_{n_{k}})-g(u)|_{s}^{s}=\int\limits_{\Omega}|g(u_{n_{k}})-g(u)|^{s}dx\stackrel{k}{\ri} %%@
0\,.$$
This end the proof of the lemma.\qed

The application $\varphi\mapsto g(x,\varphi(x))$ is the Nemitski operator
of the function $g$.

\begin{prop} Let $g:\Omega\times\RR\ri\RR$ be a Carath\'eodory function such that
$|g(x,s)|\le a+b|s|^{r-1}$ for all $(x,s)\in\Omega\times\RR$, 
with $1\le r<2N/(N-2)$ if $N\geq 2$ or $1\le r<\infty$ if $1\leq N\leq 2$. Denote 
$G(x,t)=\int_0^tg(x,s)ds$.
Let $I:H^1_0(\Omega)\ri{\mathbb R}$ be the functional characterized by
$$I(u)=\frac{1}{2}\int\limits_{\Omega}|\nabla u|^{2}dx-
\frac{\lambda_{1}}{2}\int\limits_{\Omega}u^{2}dx-\int\limits_{\Omega}G(x,u(x))dx\,.$$
Assume that
$(u_{n})_{n}\subset H^1_0(\Omega)$ has a bounded subsequence and
$I^{'}(u_{n})\ri 0$ as $n\ri\infty$. Then $(u_{n})_{n}$ has a convergent subsequence.
\end{prop}

\proof
We have
$$\langle I^{'}(u),v\rangle=\int\limits_{\Omega}\nabla u\nabla vdx-
\lambda_{1}\int\limits_{\Omega}uvdx-\int\limits_{\Omega}g(x,u(x))v(x)dx\,.$$
Denote by
$$\langle a(u),v\rangle=\int\limits_{\Omega}\nabla u\nabla vdx\,;$$
$$J(u)=\frac{\lambda_{1}}{2}\int\limits_{\Omega}u^{2}dx+\int\limits_{\Omega}G(x,u(x))dx\,.$$
It follows that
$$\langle %%@
J^{'}(u),v\rangle=\lambda_{1}\int\limits_{\Omega}uvdx+\int\limits_{\Omega}g(x,u(x))v(x)dx$$
and $ I^{'}(u)=a(u)-J^{'}(u)$. We prove that $a$ is an isomorphism
from $H^1_0(\Omega)$ onto $a(H^1_0(\Omega))$ and $J^{'}$ is
a compact operator. This assumption yields
$$u_{n}=a^{-1}\langle %%@
(I^{'}(u_{n}\rangle)+J^{'}(u_{n}))\ri\lim\limits_{n\ri\infty}a^{-1}\langle(J^{'}(u_{n})
\rangle)\,.$$
But $J^{'}$ is a compact operator and $(u_{n})_{n}$ is a bounded sequence.
This implies that $(J^{'}(u_{n}))_{n}$ has a convergent subsequence and, consequently,
$(u_{n})_{n}$ has a convergent subsequence. Assume, up to a subsequence, that
 $(u_{n})_{n}\subset H^1_0(\Omega)$ is bounded.
From the compact embedding $H^1_0(\Omega)\hookrightarrow L^{r}(\Omega)$, we can assume, passing %%@
again at a subsequence, that
 $u_{n}\ri u$ in $L^{r}(\Omega)$. We have 
\begin{eqnarray}\label{et19}
\begin{array}{ll}
\di |\!|J^{'}(u_{n})-J^{'}(u)|\!|&\di\le
\sup\limits_{|\!|v|\!|\le 1}
\left|\int\limits_{\Omega}\left(g(x,u_{n}(x))-g(x,u(x))\right)v(x)dx\right|
+\sup\limits_{|\!|v|\!|\le 1}\lambda_{1}\left|
\int\limits_{\Omega}(u_{n}-u)vdx\right|\\
&\di\le\sup\limits_{|\!|v|\!|\le 1}\int\limits_{\Omega}|g(x,u_{n}(x))-
g(x,u(x))||v(x)|dx
+\lambda_{1}\sup\limits_{|\!|v|\!|\le 1}
\int\limits_{\Omega}|(u_{n}-u)v|dx\\
&\di\leq\sup\limits_{|\!|v|\!|\le %%@
1}\left(\int\limits_{\Omega}|g(x,u_{n})-g(x,u)|^{\frac{r}{r-1}}dx\right)^{\frac{r-1}{r}}|v|_{r}\\
&\di +\lambda_{1}\sup\limits_{|\!|v|\!|\le 1}\left(\int\limits_{\Omega}|
u_{n}-u|^{2}dx\right)^{1/2}|v|_{2}\\
&\di\leq c\sup\limits_{|\!|v|\!|\le 1}
 \left(\int\limits_{\Omega}|g(x,u_{n})-g(x,u)|^{\frac{r}{r-1}}dx\right)^{\frac{r-1}{r}}|
 \!|v|\!|
+\lambda_{1}|u_{n}-u|_{2}\,.
\end{array}
\end{eqnarray}

 By Lemma \ref{brezislema} we obtain $g\in C(L^{r},L^{r/(r-1)})$. Next, since
$u_{n}\ri u$ in $L^{r}$ and $u_{n}\ri u$ in $L^{2}$, it follows by \eq{et19}
that $J^{'}(u_{n})\ri J^{'}(u)$ as $n\ri\infty$, that is, $J^{'}$ is a compact operator.
This completes our proof. \qed

Let $V$ denote the linear space spanned by $\varphi_1$. Set
$$\Gamma :=\{\gamma\in C(B,H^1_0(\Omega));\ \gamma(v)=v\,,\ \mbox{for all } v\in\partial B\} 
\qquad\mbox{and}\qquad B=\{v\in V;\ |\!|v|\!|\le R\}\,.$$
The following result has been used in the proof of Lemma \ref{deformm}.

\begin{prop}
 We have $\gamma(B)\bigcap W\not=\emptyset$, for all $\gamma\in\Gamma$.
\end{prop}

\proof
Let  $P:H^1_0(\Omega)\ri V$ be the
projection of $H^1_0$ in $V$. Then $P$ is a linear and
continuous operator. If $v\in\partial B$ then $(P\circ\gamma)(v)=P(\gamma(v))=P(v)=v$
and, consequently, $P\circ\gamma=Id$ on $\partial B$. We have
$P\circ\gamma\ ,Id\in C(B,H^1_0)$ and 
$0\not\in Id(\partial B)=\partial B$.
Using a property of the Brouwer topological degree  we obtain
$\mbox{deg}\,(P\circ\gamma,{\rm Int} B,0)=\mbox{deg}\,(Id,{\rm Int}B,0)$.
But $0\in{\rm Int}\,B$ and it follows that $\mbox{deg}\,(Id,{\rm Int}\,B,0)=1\not=0$.
So, by the existence property of the Brouwer degree, there exists $v\in{\rm Int}\,B$ such that %%@
$(P\circ\gamma)(v)=0$, that is, $P(\gamma(v))=0$.
Therefore $\gamma(v)\in W$ and this shows that $\gamma(B)\cap W\not=\emptyset$. \qed

\end{document}